\documentclass[a4paper, reqno]{amsart}

\title[{On Miura maps for $\mathcal{W}$-superalgebras}]{On Miura maps for $\mathcal{W}$-superalgebras}

\author{Shigenori Nakatsuka}
\address{Graduate School of Mathematical Sciences, The University of Tokyo, 3-8-1 Komaba, Tokyo, Japan 153-8914}
\email{nakatuka@ms.u-tokyo.ac.jp}
\address{Kavli Institute for the Physics and Mathematics of the Universe (WPI), The University of Tokyo Institutes for Advanced Study, The University of Tokyo, Kashiwa, Chiba 277-8583, Japan}
\email{shigenori.nakatsuka@ipmu.jp}

\usepackage{latexsym}
\usepackage{amsmath}
\usepackage{amsfonts}
\usepackage{amssymb}
\usepackage{amsxtra}
\usepackage{mathrsfs}
\usepackage[all]{xy}
\usepackage{color}

\newtheorem{definition}{Definition}[section]
\newtheorem{proposition}[definition]{Proposition}
\newtheorem{theorem}[definition]{Theorem}
\newtheorem{corollary}[definition]{Corollary}
\newtheorem{lemma}[definition]{Lemma}
\newtheorem{remark}[definition]{Remark}

\numberwithin{equation}{section}

\begin{document}

\maketitle

\begin{abstract}
We prove the injectivity of the Miura maps for $\mathcal{W}$-superalgebras and the isomorphisms between the Poisson vertex superalgebras obtained as the associated graded of the $\mathcal{W}$-superalgabras in terms of Li's filtration and the level 0 Poisson vertex superalgebras associated with the arc spaces of the corresponding Slodowy slices in full generality.
\end{abstract}

\section*{Introduction}
Let $\mathfrak{g}$ be a simple finite-dimensional Lie superalgebra over the field of complex numbers $\mathbb{C}$ with a non-degenerate even supersymmetric invariant bilinear form. Let $k\in \mathbb{C}$ be a complex number, $f$ a nonzero even nilpotent element of $\mathfrak{g}$ and $\Gamma\colon \mathfrak{g}=\bigoplus_{j\in\frac{1}{2} \mathbb{Z}}\mathfrak{g}_j$ a good grading for $f$. Then a vertex superalgebra $\mathcal{W}^k(\mathfrak{g},f;\Gamma)$, called the (affine) $\mathcal{W}$-superalgebra, is defined as the quantum Drinfeld-Sokolov reduction of the universal affine vertex superalgebra $V^k(\mathfrak{g})$ of $\mathfrak{g}$ at level $k$ \cite{FF0,KRW}.
Among all the even nilpotent elements, there is a distinguished one, called a principal nilpotent element (unique up to conjugation), which is characterized by the property that the dimension of its centralizer in the even subalgebra has minimal dimension. The $\mathcal{W}$-superalgebras associated with principal nilpotent elements are called the principal $\mathcal{W}$-superalgebras and in the non-super case, they enjoy the Feigin-Frenkel duality \cite{FF}:
$$\mathcal{W}^k(\mathfrak{g})\simeq \mathcal{W}^{\hspace{0.1mm}^Lk}(\hspace{0.1mm}^L \mathfrak{g}),$$
where $\hspace{0.1mm}^L \mathfrak{g}$ denotes the Langlands dual of $\mathfrak{g}$ and the levels are related by $$r(k+h^\vee)(\hspace{0.1mm}^Lk+\hspace{0.1mm}^Lh^\vee)=1.$$
Here $h^\vee$ (resp.\ $\hspace{0mm}^Lh^\vee)$ is the dual Coxeter numbere of $\mathfrak{g}$ (resp.\ $\hspace{0mm}^L \mathfrak{g})$ and $r$ is the (common) lacing number of $\mathfrak{g}$ and $\hspace{0.1mm}^L \mathfrak{g}$. Moreover, if $r=1$, i.e., $\mathfrak{g}$ is simply-laced, then there is another realization of $\mathcal{W}^k(\mathfrak{g})$:   
$$\mathcal{W}^k(\mathfrak{g})\simeq \mathrm{Com}\left(V^{k'+1}(\mathfrak{g}),V^{k'}(\mathfrak{g})\otimes L_1(\mathfrak{g})\right),$$
which is a higher-rank analogue \cite{ACL} of the Goddard-Kent-Olive construction of the Virasoro vertex algebra \cite{GKO}.
Here $L_1(\mathfrak{g})$ is the simple quotient of $V^1(\mathfrak{g})$, $\mathrm{Com}(W,V)$ denotes the coset vertex algebra of $V$ by $W$, and the levels are related by 
$$\frac{1}{k+h^\vee}+\frac{-1}{k'+h^\vee}=1.$$
Therefore, the principal $\mathcal{W}$-algebra associated with a simly-laced Lie algebra enjoys a \emph{triality}.

Recently, Gaiotto-Rap\v{c}\'{a}k \cite{GR} and Proch\'{a}zka-Rap\v{c}\'{a}k \cite{PR} proposed a generalization of the triality for $\mathcal{W}$-superalgebras and their cosets by affine vertex subalgebras in terms of 4 dimensional $\mathcal N=4$ topologically twisted super Yang-Mills theories. Since the Kazama-Suzuki coset construction of the $\mathcal N=2$ superconformal algebra appears as a very special case, these isomorphisms are expected to be efficient for the study of the representation theory of vertex superalgebras, see \cite{CGN,CL} for a mathematical approach.

In the proofs of the above triality, a vertex superalgebra homomorphism, called the {\em Miura map} \cite{A,FF2,G}
\begin{align}\label{Intro eq: Miura map} 
\mu\colon \mathcal{W}^k(\mathfrak{g},f;\Gamma)\rightarrow V^{\tau_k}(\mathfrak{g}_0)\otimes \Phi(\mathfrak{g}_{\frac12})
\end{align}
plays a fundamental role and its injectivity is very important. (See Section \ref{sec: W-superalgebra and Miura map} for details.)
The situation is expected to be the same for the generalized triality in the sense of \cite{GR,PR}. 
In general the injectivity of the map $\mu$ follows from that of its associated graded in terms of Li's filtration, see Corollary \ref{corollary: implication of injectivity}. Thus the injectivity of $\mu$ when $\mathfrak{g}$ is non-super is an immediate corollary of \cite[Corollary 8.2]{DKV}, where all the structures of the associated graded, which are called the classical $\mathcal{W}$-algebras, are determined including the explicit form of strong generators. It is natural to expect that their overall treatment of classical $\mathcal{W}$-algebras is generalized to the super setting. There is another approach \cite{F,A2}, which uses geometry behind $\mu$. The injectivity of $\mu$ is proved there when $\mathfrak{g}$ is non-super and $f$ is principal, but the proof also applies for an arbitrary $f$. For the super case, it is proven in \cite{G} by using a cohomological method when $\mathfrak{g}$ and $f$ are arbitrary but the level $k$ is generic. 
Unfortunately, it is not enough for the application to the representation theory of vertex superalgebras based on the triality since the representation theory is rich only for {\em special levels} out of generic levels.
The main theme of the paper is to improve the situation and, especially, to prove the injectivity in full generality. 

In Section \ref{sec: Formal supergroups}, we extend some results on unipotent algebraic groups to the super settings in view of the study of $\mathcal{W}$-superalgebras. We explain the correspondence between formal supergroups and finite dimensional Lie superalgebras \cite{Ma1,Se} by using the Campbell-Hausdorff formal supergroups. This correspondence gives the one between unipotent algebraic supergroups and finite dimensional nilpotent Lie superalgebras \cite{MO}. This implies that every unipotent algebraic supergroup is an affine superspace (i.e., $\mathbb{C}^{p|q}$ for some $p,q\in \mathbb{Z}_{\geq0}$) as is well-known in the non-super cases, cf \cite{Mil}. Let $P$ be an unipotent algebraic supergroup and  $\mathfrak{p}$ the corresponding Lie superalgebra. Then we have a left $\mathfrak{p}$-module structure on the coordinate ring $\mathbb{C}[P]$. We prove that under a certain condition, the Lie superalgebra cohomology $H^n(\mathfrak{p};\mathbb{C}[P])$ is isomorphic to $\delta_{n,0}\mathbb{C}$. The proof is very standard: we use a filtration on the complex, which gives a spectral sequence $E_1\simeq H^\bullet_{\mathrm{dR}}(P)\Rightarrow H^\bullet(\mathfrak{p};\mathbb{C}[P])$ and then use Poincar\'{e} lemma $H^n_{\mathrm{dR}}(P)\simeq \delta_{n,0}$ since $P$ is an affine superspace. 

In Section \ref{sec: W-superalgebra and Miura map}, after reviewing the definition of $\mathcal{W}$-superalgebras and the construction of the Miura maps \eqref{Intro eq: Miura map}, we extend some results of the principal $\mathcal{W}$-algebras to the super cases, following \cite{A}. 
In this section, the readers are supposed to be familiar with vertex superalgebras (see e.g.\ \cite{Ka2}). Since we will also use Poisson vertex superalgebras, we will use the language of $\lambda$-bracket, see e.g., \cite{Ka3} for details. 
We first consider Li's filtrations \cite{L} on $\mathcal{W}^k(\mathfrak{g},f;\Gamma)$ and $V^{\tau_k}(\mathfrak{g}_0)\otimes \Phi(\mathfrak{g}_{\frac12})$, respectively. The associated graded vector superspaces admit a natural structure of Poisson vertex superalgebras (PVAs) and $\mu$ induces a PVA homomorphism
$$\bar{\mu}\colon \operatorname{gr}_F\mathcal{W}^k(\mathfrak{g},f;\Gamma)\rightarrow \operatorname{gr}_F\left(V^{\tau_k}(\mathfrak{g}_0)\otimes \Phi(\mathfrak{g}_{\frac12})\right).$$
The injectivity of $\mu$ reduces to that of $\bar{\mu}$ (Lemma \ref{injectivity: general setting}).
Let $\{e,h,f\}\subset \mathfrak{g}$ be an $\mathfrak{sl}_2$-triple in the even subalgebra of $\mathfrak{g}$ containing the nilpotent element $f$. We write $\mathfrak{g}_{\geq-1/2}=\oplus_{j\geq-1/2}\mathfrak{g}_j$, $\mathfrak{g}_+=\oplus_{j>0}\mathfrak{g}_j$, and denote by $G_+$ the unipotent algebraic supergroup whose Lie superalgebra is $\mathfrak{g}_+$. Then the affine subspace $f+\mathfrak{g}_{\geq-1/2}$ is stable under the adjoint action of $G_+$ and, moreover, admits a $G_+$-equivariant isomorphism of affine supervarieties
$$\mathcal{S}_f\times G_+\simeq f+\mathfrak{g}_{\geq-1/2}, \quad (X,g)\mapsto g^{-1}Xg,$$
where $\mathcal{S}_f=f+\mathfrak{g}^e\subset \mathfrak{g}$, called the {\em Slodowy slice} of $f$ (Proposition \ref{slodowy slice}). This is a generalization to the super setting of \cite[eq.(4)]{A2} and \cite[Lemma 2.1]{GG}.
We have a Poisson structure on $\mathbb{C}[f+\mathfrak{g}_{\geq-1/2}]$ (see Section \ref{sec: Arc space of Slodowy slice}), which restricts to $\mathbb{C}[\mathcal{S}_f]\simeq \mathbb{C}[f+\mathfrak{g}_{\geq-1/2}]^{G_+}$. If $\Gamma$ is a $\mathbb{Z}$-grading, then the Poisson supervariety $\mathcal{S}_f$ is also obtained as the Hamiltonian reduction of $\mathfrak{g}$ equipped with the Kostant-Kirillov Poisson structure with respect to $G_+$ and character $\chi=(f|\text{-})$. It turns out that we have an isomorphism $\operatorname{gr}_F\mathcal{W}^k(\mathfrak{g},f;\Gamma)\simeq \mathbb{C}[J\mathcal{S}_f]$ as PVAs (Proposition \ref{proposition: classical W-superalgebra}) where $J\mathcal{S}_f$ is the arc space of $\mathcal{S}_f$ (see Section \ref{Superscheme of formal arcs}) and $\mathbb{C}[J\mathcal{S}_f]$ is equipped with the (level 0) PVA structure induced by the Poisson structure on $\mathbb{C}[\mathcal{S}_f]$. On the other hand, we identify $\operatorname{gr}_F\left(V^{\tau_k}(\mathfrak{g}_0)\otimes \Phi(\mathfrak{g}_{\frac12})\right)\simeq \mathbb{C}[J(f+\mathfrak{g}_{-1/2}\oplus \mathfrak{g}_0)]$ and thus
$$\bar{\mu}\colon \mathbb{C}[J\mathcal{S}_f]\rightarrow \mathbb{C}[J(f+\mathfrak{g}_{-1/2}\oplus \mathfrak{g}_0)].$$
This map is induced by its finite dimensional analogue 
$$\bar{\mu}_{\mathrm{fin}}\colon \mathbb{C}[\mathcal{S}_f]\simeq \mathbb{C}[f+\mathfrak{g}_{\geq-1/2}]^{G_+}\rightarrow \mathbb{C}[f+\mathfrak{g}_{-1/2}\oplus \mathfrak{g}_0],$$
which is just the restriction map. We prove the injectivity of $\bar{\mu}_{\mathrm{fin}}$ (Proposition \ref{classical finite Miura}), which immediately implies the injectivity of $\bar{\mu}$ and thus of the original Miura map $\mu$ (Theorem \ref{th: injectivity of Miura}).

\vspace{3mm}

{\it Acknowledgments}\quad 
The author would like to express his gratitude to Thomas Creutzig and Naoki Genra for suggesting the problem and useful discussions, to his supervisor Masahito Yamazaki for reading the draft, to Yuto Moriwaki for useful discussions on algebraic supergroups and to Daniel Valeri for informing him of his result with A. De Sole and V. Kac of the injectivity of the Miura map in the purely even cases.

\section{Unipotent algebraic supergroups}\label{sec: Formal supergroups}
\subsection{Formal supergroups}\label{sec: fromal supergroups}
We review the contravariant equivalence between formal supergroups and Lie superalgebras, following \cite{Se,Ma1}. The results in Section \ref{sec: fromal supergroups}-\ref{sect: Unip-Nilp} are proved in much more generality in \cite{Ma1}.
Given a commutative $\mathbb{C}$-superalgebra $A$, we denote by $A=A_{\bar{0}}\oplus A_{\bar{1}}$ the parity decomposition and by $\bar{a}$ the parity of $a\in A$. In particular, given a parity-homogeneous basis $\{x_\alpha\}_{\alpha\in S}$ of $A$ with an index set $S$, the parity $\bar{x}_\alpha$ is also denoted by $\bar{\alpha}$. The multiplication $m\colon A\times A\rightarrow A$ is also denoted by $ab=m(a,b)$.

Let $\{x_\alpha\}_{\alpha\in S}$ be a set equipped with parity, i.e., $S=S_{\bar{0}}\sqcup S_{\bar{1}}$ and $\bar{x}_{\alpha}=\bar{\alpha}$.
The ring of (commutative) superpolynomials in the variables $\{x_\alpha\}_{\alpha\in S}$ is the commutative $\mathbb{C}$-superalgebra
$$\mathbb{C}[x_\alpha\mid \alpha\in S]:=\mathbb{C}[x_\alpha\mid \alpha\in S_{\bar{0}}]\otimes \bigwedge \ _{x_\beta, (\beta\in S_{\bar{1}})},$$
which is the tensor product of the polynomial ring generated by the variables $\{x_\alpha\}_{\alpha\in S_{\bar{0}}}$ and the exterior algebra generated by the variables $\{x_\alpha\}_{\alpha\in S_{\bar{1}}}$.

We denote by $\hat{R}_{p|q}$ a formal power series ring 
$$\hat{R}_{p|q}=\mathbb{C} [\![x_1,\cdots,x_p]\!]\otimes \bigwedge \ _{\phi_1,\cdots,\phi_q}$$
with $p$ even variables and $q$ odd variables, and by $\hat{\mathfrak{m}}=(x_1,\cdots,x_p,\phi_1,\cdots,\phi_q)\subset \hat{R}_{p|q}$ its unique maximal ideal. It is complete in the linear topology whose basis of open neighborhoods around 0 is $\{\hat{\mathfrak{m}}^n\}_{n=0}^\infty$. Let $m\colon \hat{R}_{p|q}\times \hat{R}_{p|q}\rightarrow \hat{R}_{p|q}$ denote the product and $u\colon \mathbb{C}\rightarrow \hat{R}_{p|q}$ the unit morphism.
We denote by $\hat{R}_{p|q}\widehat{\otimes}\hat{R}_{p|q}$ the completed tensor product
$\hat{R}_{p|q}\widehat{\otimes}\hat{R}_{p|q}:=\varprojlim_{n,m}\hat{R}_{p|q}/\hat{\mathfrak{m}}^n\otimes \hat{R}_{p|q}/\hat{\mathfrak{m}}^m$. Then we have a natural isomorphism
$$\hat{R}_{p|q}\widehat{\otimes}\hat{R}_{p|q}\simeq \mathbb{C}\left[\!\!\left[x_1^{(1)},\cdots,x_p^{(1)},x_1^{(2)},\cdots,x_p^{(2)}\right]\!\!\right]\otimes \bigwedge \ _{\phi_1^{(1)},\cdots,\phi_q^{(1)},\phi_1^{(2)},\cdots,\phi_q^{(2)}}.$$

A {\it formal supergroup} is a super bialgebra $(\hat{R}_{p|q},\Delta,\epsilon)$ with continuous homomorphisms
\begin{itemize}
\item
coproduct  \ $\Delta\colon \hat{R}_{p|q}\rightarrow \hat{R}_{p|q}\widehat{\otimes}\hat{R}_{p|q}$,
\item
counit \ $\epsilon\colon \hat{R}_{p|q}\rightarrow \mathbb{C}$.
\end{itemize}
It has an unique continuous homomorphism $S\colon \hat{R}_{p|q}\rightarrow \hat{R}_{p|q}$ such that 
$$m\circ(S\otimes 1)\circ \Delta=\epsilon \circ u=m\circ(1\otimes S)\circ \Delta,$$
called the antipode. Thus formal supergroups are Hopf superalgebras, see \cite{CCF}. 
\begin{remark}\label{taking germ at e}
Let $G$ be an algebraic supergroup with the identity $e$. Then the germ of the structure sheaf $\mathcal{O}_G$ at $e$, which we denote by $\mathcal{O}_{G,e}$, is a local ring. Let $\mathfrak{m}_{G,e}$ denote its unique maximal ideal. Then the completion
$$\widehat{\mathcal{O}}_{G,e}:=\varprojlim \mathcal{O}_{G,e}/\mathfrak{m}_{G,e}^n$$
has a natural structure of formal supergroup.
\end{remark}

A morphism of formal supegroups $(\hat{R}_{p|q},\Delta,\epsilon)$ and $(\hat{R}_{p'|q'},\Delta',\epsilon')$
is a morphism of Hopf superalgebras 
$$F\colon \hat{R}_{p|q}\rightarrow \hat{R}_{p'|q'},$$
which is continuous in the linear topology. We denote by $\underbar{FG}$ the category of formal supergroups.

\subsection{Point distributions of formal supergroups}
Let $F=(\hat{R}_{p|q},\Delta,\epsilon)$ be a formal supergroup and $\hat{\mathfrak{m}}\subset \hat{R}_{p|q}$ its maximal ideal.
A {\it point distribution} is a continuous linear map
$$\xi\colon \hat{R}_{p|q}\rightarrow \mathbb{C},$$
where $\mathbb{C}$ is equipped with the discrete topology. Let $U_F:=\mathrm{Hom}_\mathbb{C}^{c.}(\hat{R}_{p|q},\mathbb{C})$ denote the vector superspace of point distributions. Since any point distribution factors 
$$\hat{R}_{p|q}\twoheadrightarrow \hat{R}_{p|q}/\hat{\mathfrak{m}}^\ell\rightarrow \mathbb{C}$$
for some $\ell\in \mathbb{Z}_{\geq0}$, we have
$U_F=\varinjlim_{\ell}\mathrm{Hom}_\mathbb{C}(\hat{R}_{p|q}/\hat{\mathfrak{m}}^\ell,\mathbb{C})$. It has an associative algebra structure by the dual maps 
$$\Delta^*\colon U_F\otimes U_F\rightarrow U_F,\quad \epsilon^*\colon \mathbb{C}\rightarrow U_F$$
of $\Delta$ and $\epsilon$, respectively. 
Note that these maps are well-defined since $\Delta$ and $\epsilon$ are continuous. Moreover, $U_F$ has a Hopf superalgebra structure $F^\vee=(U_F, m^*,u^*,S^*)$.

Let us describe $F^\vee$ more explicitly. 
The elements $X^\alpha=x_1^{i_1}\cdots x_p^{i_p}\phi_1^{j_1}\dots \phi_q^{j_q}\in \hat{R}_{p|q}$, ($\alpha=(\alpha_1,\alpha_2)=(i_1,\cdots,i_p,j_1,\cdots,j_q)\in S=\mathbb{Z}_{\geq0}^p\times \{0,1\}^q$), form a topological basis of $\hat{R}_{p|q}$. We denote by $\xi_\alpha$, ($\alpha\in S$), its dual basis of $U_F$, i.e., $\xi_\alpha$ satisfies $\xi_\alpha(X^\beta)=\delta_{\alpha,\beta}$. 
Then $\xi_\alpha\otimes \xi_\beta$, ($\alpha,\beta\in S$), form a basis of $U_F\otimes U_F$. Here $\xi_\alpha\otimes \xi_\beta$ satisfies 
$$\xi_\alpha\otimes \xi_\beta(X^{\alpha'}\otimes X^{\beta'})=(-1)^{\bar{\alpha'}\bar{\beta}}\xi_\alpha(X^{\alpha'})\xi_{\beta}(X^{\beta'}).$$
 In the following, we extend the basis $\xi^\alpha\in U_F$, ($\alpha\in S=\mathbb{Z}_{\geq0}^p\times \{0,1\}^q$), 
 to $\xi\in U_F$, ($\alpha\in \widetilde{S}:=\mathbb{Z}_{\geq0}^p\times \mathbb{Z}_{\geq0}^q$), by setting 0 for $\alpha\in \widetilde{S}\backslash S$. We define a degree on $U_F$ by $\mathrm{deg}(\xi_\alpha)=|\alpha|:=|\alpha_1|+|\alpha_2|$.

Then we have
\begin{align}\label{counit and unit on UF}
\begin{array}{cccccccc}
\epsilon^*\colon& \mathbb{C}&\rightarrow& U_F,&\quad u^*\colon &U_F&\rightarrow &\mathbb{C}\\
& 1&\mapsto &\xi_0,&\quad &\xi_\alpha&\mapsto&\delta_{\alpha,0}.
\end{array}
\end{align}
The product $F^*\colon U_F\times U_F\rightarrow U_F$ satisfies  
\begin{align}\label{product on UF}
\xi_\alpha\cdot \xi_\beta=(-1)^{\bar{\alpha}\bar{\beta}}
\binom{\alpha_1+\beta_1}{\alpha_1}
\xi_{\alpha+\beta}+\text{lower degree terms},
\end{align}
where
$$\binom{p}{q}:=\prod_{i=1}^\ell\binom{p_i}{q_i},\quad p=(p_1,\cdots,p_\ell),\ q=(q_1,\cdots,q_\ell),$$
and the lower degree terms do not contain constant term $\mathbb{C} \xi_0$.  
In particular, it follows that the elements $\xi_{\alpha}$ with $|\alpha|=1$ generate $U_F$ as a superalgebra.
The coproduct $m^*$ is the unique superalgebra homomorphism such that 
\begin{align}\label{coproduct on UF}
m^*(\xi_\alpha)=\xi_\alpha\otimes 1+1\otimes \xi_\alpha
\end{align}
for $\alpha\in S$ such that $|\alpha|=1$.
\begin{lemma}[{\cite[Lemma 4.1]{Ma1}}]
The bilinear map 
$$[\cdot,\cdot]\colon U_F\times U_F\rightarrow U_F,\quad (a,b)\mapsto ab-(-1)^{\bar{a}\bar{b}}ba$$
defines a Lie superalgebra structure on the subspace $U_F^1=\oplus_{|\alpha|=1}\mathbb{C} \xi_\alpha$.
\end{lemma}
\proof
Since $U_F$ is an associative superalgebra, it suffices to prove that $[\cdot,\cdot]$ preserves the subspace $U_F^1$. But it follows from 
 \eqref{product on UF}.
\endproof
Let $\mathcal{U}(U_F^1)$ denote the universal enveloping superalgebra of the Lie superalgebra $U_F^1$. Recall that it admits a Hopf superalgebra structure with coproduct
$$\mathcal{U}(U_F^1)\rightarrow \mathcal{U}(U_F^1)\otimes \mathcal{U}(U_F^1),\quad X\mapsto X \otimes 1+1\otimes X,\quad (X\in U_F^1),$$
counit 
$$\mathcal{U}(U_F^1)\twoheadrightarrow \mathcal{U}(U_F^1)/(U_F^1)\simeq \mathbb{C},$$
and antipode
$$\mathcal{U}(U_F^1)\rightarrow \mathcal{U}(U_F^1),\quad X\mapsto -X,\quad (X\in U_F^1).$$
Here $(U_F^1)\subset \mathcal{U}(U_F^1)$ denotes the left ideal generated by the subspace $U_F^1$.
\begin{proposition}[{\cite[Proposition 4.8]{Ma1}}]\label{formal supergroup to Lie superalgebra}
The Lie superalgebra homomorphism $U_F^1\hookrightarrow U_F$ induces a $\mathbb{C}$-superalgebra homomorphism $\mathcal{U}(U_F^1)\rightarrow U_F$. It is an isomorphism of Hopf superalgebras.
\end{proposition}
\proof
It is obvious that we have a  $\mathbb{C}$-superalgebra homomorphism $\mathcal{U}(U_F^1)\rightarrow U_F$. Let $\iota$ denote this map. 
Then it follows from \eqref{counit and unit on UF} and \eqref{coproduct on UF} that $\iota$ is a super bialgebra homomorphism. Since an antipode on a super bialgebra is unique if it exists, $\iota$ intertwines the antipodes of $\mathcal{U}(U_F^1)$ and $U_F$. Thus $\iota$ is a homomorphism of Hopf superalgebras. Since the elements $\xi_\alpha$with $\alpha\in S$ such that $|\alpha|=1$ generate $U_F$ as a $\mathbb{C}$-superalgebra, $\iota$ is surjective. Now the PBW type theorem for $\mathcal{U}(U_F^1)$ implies that $\iota$ is an isomorphism.
\endproof
\subsection{Equivalence between formal supergroups and Lie superalgebras}
Let $\underbar{LA}$ denote the category of finite dimensional Lie superalgebras. Taking the universal enveloping superalgebra $\mathfrak{a}\mapsto \mathcal{U}(\mathfrak{a})$ for a Lie superalgebra $\mathfrak{a}$ defines a functor
$$\mathcal{U}\colon \underbar{LA}\rightarrow \underbar{BA}$$
from $\underbar{LA}$ to the category of super bialgebras $\underbar{BA}$. 
The functor $\mathcal{U}$ is fully faithful and its quasi-inverse is given by taking the primitive elements:
$$\mathcal{P}\colon \underbar{BA}\rightarrow \underbar{LA},\quad (R,\Delta,\epsilon)\mapsto \mathcal{P}(R):=\{a\in R\mid \Delta(a)=a\otimes 1+1\otimes a\}.$$
The Lie superalgebra structure is given by $(a,b)\mapsto ab-(-1)^{\bar{a}\bar{b}}ba$ for $a,b\in \mathcal{P}(R)$
Thus $\underbar{LA}$ is regarded as a full subcategory of $\underbar{BA}$. 
By Proposition \ref{formal supergroup to Lie superalgebra}, the association $F\mapsto F^\vee$ for a formal supergroup $F$ gives a contravariant functor
$(\text{-})^\vee\colon \underbar{FG}\rightarrow \underbar{LA}$
from the category of formal supergroups introduced at the end of Section 1.1, to the category of finite dimensional Lie superalgebras.
\begin{theorem}[{\cite[Theorem 4.4]{Ma1}}]
The functor $(\text{-})^\vee\colon \underbar{FG}\rightarrow \underbar{LA}$ gives a contravariant equivalence of categories.
\end{theorem}
\noindent A quasi-inverse is given again by taking the dual vector superspace
\begin{align}\label{equiv from LA to FG}
(\text{-})^\vee\colon \underbar{LA}\rightarrow \underbar{FG},\quad \mathfrak{a}\mapsto \mathrm{Hom}_{\mathbb{C}}(\mathcal{U}(\mathfrak{a}),\mathbb{C}).
\end{align}
The super coalgebra (algebra) structure on $\mathrm{Hom}_{\mathbb{C}}(\mathcal{U}(\mathfrak{a}),\mathbb{C})$ is induced by the super algebra (coalgebra) structure on $\mathrm{Hom}_{\mathbb{C}}(\mathcal{U}(\mathfrak{a}),\mathbb{C})$. In the purely even setting, the equivalence between the categories of formal groups and of finite dimensional Lie algebras is a special case of {\it Cartier duality} between the category of linearly compact commutative Hopf algebras and the category of cocommutative Hopf algebras, see \cite{D} for details.

The quasi-inverse \eqref{equiv from LA to FG} is known to be naturally isomorphic to the functor of taking the Campbell-Hausdorff formal supergroup. 
By viewing a finite dimensional Lie superalgebra $\mathfrak{a}$ as an affine superscheme, we denote by $\mathbb{C}[\mathfrak{a}]$ its coordinate ring. We identify $\mathbb{C}[\mathfrak{a}]$ with the symmetric superalgebra $S(\mathfrak{a}^*)$ of the dual vector superspace $\mathfrak{a}^*$ of $\mathfrak{a}$. Then the Lie superbracket $[\cdot,\cdot]$ on $\mathfrak{a}$ induces a superalgebra homomorphism
\begin{align}\label{cobracket}
D\colon S(\mathfrak{a}^*)\rightarrow S(\mathfrak{a}^*)\otimes S(\mathfrak{a}^*),
\end{align}
called the co(super)bracket of $\mathfrak{a}$.
Consider the completion $\hat{S}(\mathfrak{a}^*)=\varprojlim_{n}S(\mathfrak{a}^*)/\mathcal{I}^n$ where $\mathcal{I}$ is the argumentation ideal and the completed tensor product 
$$\hat{S}(\mathfrak{a}^*)\widehat{\otimes}\hat{S}(\mathfrak{a}^*)=\varprojlim_{n,m}S(\mathfrak{a}^*)/\mathcal{I}^n\otimes S(\mathfrak{a}^*)/\mathcal{I}^m.$$ 
Then the vector superspace $\mathrm{Hom}_{\mathbb{C}}(\mathfrak{a}^*,\hat{S}(\mathfrak{a}^*)\widehat{\otimes}\hat{S}(\mathfrak{a}^*))$ admits a Lie superalgebra structure by
\begin{align}
[f,g]:=m\circ (f\otimes g)\circ D,\quad f,g\in \mathrm{Hom}_{\mathbb{C}}\left(\mathfrak{a}^*,\hat{S}(\mathfrak{a}^*)\widehat{\otimes}\hat{S}(\mathfrak{a}^*)\right),
\end{align}
where $m$ denotes the multiplication of $\hat{S}(\mathfrak{a}^*)\widehat{\otimes} \hat{S}(\mathfrak{a}^*)$. Define a superalgebra homomorphism 
$\Delta\colon \hat{S}(\mathfrak{a}^*)\rightarrow \hat{S}(\mathfrak{a}^*)\widehat{\otimes}\hat{S}(\mathfrak{a}^*)$
by 
\begin{align}\label{coproduct of fg}
\Delta(X)=\sum_{n=1}^{\infty}\frac{(-1)^{n+1}}{n}\sum_{\begin{subarray}{c}(p_i,q_i)>0\\ i=1,\cdots,n\end{subarray}}\frac{[\iota_1^{p_1},\iota_2^{q_1},\cdots, \iota_1^{p_n},\iota_2^{q_n}](X)}{\sum (p_i+q_i)\prod p_i!q_i!},\quad (X\in \mathfrak{a}^*),
\end{align}
where the summation $(p_i,q_i)$ is over $\mathbb{Z}_{\geq0}^2\backslash \{(0,0)\}$, $\iota_i$ denotes the natural inclusion
$$\iota_i\colon \mathfrak{a}^*\rightarrow \mathfrak{a}^*\otimes \mathbb{C} \oplus \mathbb{C}\otimes \mathfrak{a}^*\subset \hat{S}(\mathfrak{a}^*)\widehat{\otimes}\hat{S}(\mathfrak{a}^*),$$
for $i=1,2$ and 
\begin{align*}
[\iota_1^{p_1},\iota_2^{q_1},&\cdots, \iota_1^{p_n},\iota_2^{q_n}]\\
&=\begin{cases}
\underbrace{[\iota_1,\cdots,[\iota_1}_{p_1},\underbrace{[\iota_2,\cdots,[\iota_2}_{q_1}\cdots
\underbrace{[\iota_1,\cdots,[\iota_1}_{p_n},
\underbrace{[\iota_2,\cdots,[\iota_2}_{q_n-1},\iota_2],
& (q_n\neq0),\\
\underbrace{[\iota_1,\cdots,[\iota_1}_{p_1},\underbrace{[\iota_2,\cdots,[\iota_2}_{q_1}\cdots
\underbrace{[\iota_1,\cdots,[\iota_1}_{p_n-1},\iota_1],
& (q_n=0).
\end{cases}
\end{align*}
Then by \cite[Proposition 4.11]{Ma1}, $(\hat{S}(\mathfrak{a}^*),\Delta,\epsilon)$ defines a formal supergroup where $\epsilon\colon \hat{S}(\mathfrak{a}^*)\rightarrow \mathbb{C}$ is the obvious counit. Moreover, by the proof of \cite[Theorem 4.4]{Ma1}, the functor 
$$CH\colon \underbar{LA}\rightarrow \underbar{FG},\quad \mathfrak{a}\mapsto (\hat{S}(\mathfrak{a}^*),\Delta,\epsilon)$$
is a quasi-inverse of $(\text{-})^\vee\colon \underbar{FG}\rightarrow \underbar{LA}$. 

\subsection{Unipotent algebraic supergroups and nilpotent Lie superalgebras}\label{sect: Unip-Nilp}

An affine algebraic supergroup $P$ is called {\it unipotent} if the following equivalent conditions hold:
\begin{itemize}
\item
 the coordinate ring $\mathbb{C}[P]$ is irreducible as a coalgebra, i.e., $\mathbb{C}[P]$ has a unique simple subcoalgebra, which is $\mathbb{C}$,
\item 
the isomorphism classes of simple rational $P$-supermodules are trivial modules $\{\mathbb{C}^{1|0},\mathbb{C}^{0|1}\}$,
\end{itemize}
see \cite[Definition 2.9]{Ma2}. Let $\underbar{Unip-AG}$ denote the full subcategory of the category of affine algebraic supergroups consisting of unipotent affine algebraic supergroups. Then we have a functor 
$$\hat{\mathcal{O}}_{\text{-},e}\colon \underbar{Unip-AG}\rightarrow \underbar{FG},\quad P\mapsto \hat{\mathcal{O}}_{P,e},$$
by Remark \ref{taking germ at e} and 
$$T_e\text{-}\colon \underbar{Unip-AG}\rightarrow \underbar{LA},\quad P\mapsto T_eP.$$
The latter one is just taking the Lie superalgebra on the tangent space at the identity $e$, see \cite{CCF} for details. 
By construction, we have a natural isomorphism 
\begin{align}\label{tangent of unipotents}
(\text{-})^\vee\circ \hat{\mathcal{O}}_{\text{-},e}\simeq T\text{-}_e\colon \underbar{Unip-AG}\rightarrow \underbar{LA}.
\end{align}

A finite dimensional Lie superalgebra $\mathfrak{p}$ is called {\it nilpotent} if the descending central series 
\begin{align}\label{descending central series}
\mathfrak{p}^1=\mathfrak{p},\quad \mathfrak{p}^n=[\mathfrak{p},\mathfrak{p}^{n-1}],\ (n>1),
\end{align}
vanishes, i.e., $\mathfrak{p}^n=0$ for some $n$. Let $\underbar{Nil-LA}\subset \underbar{LA}$ denote the full subcategory consisting of (finite dimensional) nilpotent Lie superalgebras. 
\begin{proposition}[\cite{MO}]\label{Unip-Nilp}
The functor $T_e\text{-}\colon \underbar{Unip-AG}\rightarrow \underbar{Nil-LA}$ gives an equivalence of categories.
\end{proposition}
Let $D\colon \mathfrak{p}^*\rightarrow \mathfrak{p}^*\otimes \mathfrak{p}^*$ denote the cobracket of $\mathfrak{p}$. If $\mathfrak{p}$ is nilpotent, then the map $D^n$ defined by 
$$D^n=(D\otimes \underbrace{1\otimes\cdots \otimes1}_{n-1})\circ D^{n-1},\quad n\geq1,$$
vanishes, i.e., $D^n=0$ fore some $n\geq1$. Thus the coproduct \eqref{coproduct of fg} of the corresponding formal group $CH(\mathfrak{p})$ stabilize the supersymmetric algebra $S(\mathfrak{p}^*)$:
\begin{align}\label{coproduct of CH group in the nilpotent case}
\Delta\colon S(\mathfrak{p}^*)\rightarrow S(\mathfrak{p}^*)\otimes S(\mathfrak{p}^*).
\end{align}
Equivalently, $\mathrm{Spec}(S(\mathfrak{p}^*))$ has a structure of affine algebraic supergroup. By Proposition \ref{Unip-Nilp} and \eqref{tangent of unipotents}, $\mathrm{Spec}(S(\mathfrak{p}^*))$ is an unipotent affine algebraic supergroup and, conversely, any unipotent affine algebraic supergroup is obtained in this way. Since $\mathrm{Spec}(S(\mathfrak{p}^*))$ is an affine superspace, we obtain the following assertion. 
\begin{corollary}\label{unipotents are affine spaces}
Any unipotent affine algebraic supergroup is isomorphic to its Lie  superalgebra as affine supercpaces and thus $\mathbb{C}^{p|q}$ for some $p,q\in \mathbb{Z}_{\geq0}$.
\end{corollary}
\subsection{Lie superalgebra cohomology}
Let $P$ be an affine algebraic supergroup with Lie superalgebra $\mathfrak{p}=T_eP$.
For any commutative superalgebra $A$ and an $A$-point $g\in P(A)$, we have a right multiplication 
$$R_g\colon P(A)\rightarrow P(A),\quad h\mapsto hg.$$
This induces a left $\mathfrak{p}$-module structure on the coordinate ring
\begin{align}\label{regular rep}
R\colon \mathfrak{p}\rightarrow \mathrm{End}_\mathbb{C}(\mathbb{C}[P]).
\end{align}
We call $P$ positively graded if $\mathfrak{p}$ is $\mathbb{Z}_{>0}$-graded:
$$\mathfrak{p}=\bigoplus_{n>0}\mathfrak{p}_n,\quad [\mathfrak{p}_n,\mathfrak{p}_m]\subset \mathfrak{p}_{n+m}.$$
In this case, $\mathfrak{p}^*$ is naturally $\mathbb{Z}_{<0}$-graded by $\operatorname{wt}\colon \mathfrak{p}^*=\oplus_{n<0}(\mathfrak{p}^*)_n$ where $(\mathfrak{p}^*)_n=\operatorname{Hom}_\mathbb{C}(\mathfrak{p}_{-n},\mathbb{C})$ for $n<0$. We extend the grading on $\mathfrak{p}^*$ to $S(\mathfrak{p}^*)$ multiplicatively: $\operatorname{wt}(ab)=\operatorname{wt}(a)+\operatorname{wt}(b)$. Then the cobracket $D$ in \eqref{cobracket} preserves the grading and so does the coproduct $\Delta\colon \mathbb{C}[P]\rightarrow \mathbb{C}[P]\otimes \mathbb{C}[P]$ by \eqref{coproduct of fg}.
\begin{proposition}\label{cohomology of gra unipotent group}
For a positively graded unipotent affine algebraic supergroup $P$,
$$H^n(\mathfrak{p}; \mathbb{C}[P])\simeq \delta_{n,0}\mathbb{C}.$$
\end{proposition}
\proof
Although the assertion is well-known, we include a proof for the completeness of the paper. Since $\mathfrak{p}$ is nilpotent, we have
$$\mathfrak{p}\supset \mathfrak{p}^1=[\mathfrak{p},\mathfrak{p}]\supset  \mathfrak{p}^2=[\mathfrak{p},[\mathfrak{p},\mathfrak{p}]]\supset \cdots \supset \mathfrak{p}^m=0$$ 
for some $m\in \mathbb{Z}_{\geq0}$. We write $N=\dim \mathfrak{p}$ and take a basis $\{v_\alpha\}_{\alpha=1}^N$ of $\mathfrak{p}$ such that $\{v_i\}_{i=1}^{\mathrm{dim}\mathfrak{p}/\mathfrak{p}^1}$ gives a basis of $\mathfrak{p}/\mathfrak{p}^1$, $\{v_i\}_{i=1+\mathrm{dim}\mathfrak{p}/\mathfrak{p}^1}^{\mathrm{dim}\mathfrak{p}/\mathfrak{p}^2}$ gives a basis of $\mathfrak{p}^1/\mathfrak{p}^2$, and so on. Let $c_{\alpha \beta}^\gamma$ denote the structure constants of $\mathfrak{p}$, i.e.,
$[v_\alpha,v_\beta]=\sum_\gamma c_{\alpha \beta}^\gamma v_\gamma$. Note that $c_{\alpha \beta}^\gamma$ is non-zero only if $\gamma \geq \alpha,\beta$. The Chevalley-Eilenberg complex of $\mathfrak{p}$ with coefficients in $\mathbb{C}[P]$ is 
$$C^\bullet (\mathfrak{p};\mathbb{C}[P])=\mathbb{C}[P]\otimes S(\Pi \mathfrak{p}^*),$$
\begin{align}\label{differential}
d=\sum_{\alpha}(-1)^{\bar{\alpha}}R(v_\alpha)\otimes \varphi^\alpha-\frac{1}{2}\sum_{\alpha,\beta,\gamma}(-1)^{\bar{\alpha}\bar{\gamma}}c_{\alpha \beta}^\gamma \varphi_\gamma \varphi^\alpha \varphi^\beta.
\end{align}
Here $\Pi \mathfrak{p}^*$ denotes the parity-reversed superspace of $\mathfrak{p}^*$, spanned by $\{\varphi^\alpha\}_{\alpha=1}^N$ with parity $\bar{\varphi}^\alpha=\bar{\alpha}+1\in \mathbb{Z}_{2}$. The symbol $\varphi^\alpha$ in the differential $d$ is the multiplication by $\varphi^\alpha$ and $\varphi_\alpha$ is the contracting operator for $\varphi^\alpha$, i.e., $\varphi_\alpha(\varphi^\beta)=\delta_{\alpha,\beta}$, which has the same parity as $\varphi^\alpha$, (see e.g.\ \cite{DK}).
By \eqref{coproduct of CH group in the nilpotent case}, we may take $\mathbb{C}[\mathfrak{p}]$ as the coordinate ring $\mathbb{C}[P]$ of $P$. Let $\{x_\alpha\}_{\alpha=1}^N$ denote the linear coordinates of $\mathfrak{p}$ corresponding to the basis $\{v_\alpha\}_{\alpha=1}^N$. Then $\mathbb{C}[\mathfrak{p}]=\mathbb{C}[x_\alpha\mid \alpha=1,\cdots, N]$. By \eqref{coproduct of fg}, \eqref{regular rep} is expressed as 
$$R(v_\alpha)=\sum_{\beta}\left(\epsilon^{(2)}\frac{\partial}{\partial x_\alpha^{(2)}}\Delta(x_\beta)\right)\frac{\partial}{\partial x_\beta},$$
where 
$$\epsilon^{(2)}\frac{\partial}{\partial x_\alpha^{(2)}} \Delta(x_\beta):=\sum x_{\beta(1)}\epsilon\left(\frac{\partial}{\partial x_\alpha}x_{\beta(2)}\right)$$
and $\Delta(x_\beta)=\sum x_{\beta(1)}\otimes x_{\beta(2)}$. Thus the map $R$ preserves the grading. Extend the grading on $\mathbb{C}[P]$ to $C^\bullet(\mathfrak{p};\mathbb{C}[P])$ by 
$\operatorname{wt}(\varphi^\alpha)=\operatorname{wt}(x_\alpha)$. Then we have the grading decomposition 
\begin{align}\label{degree decomposition}
C^\bullet(\mathfrak{p};\mathbb{C}[P])=\bigoplus_{n\leq0}C^\bullet_n(\mathfrak{p};\mathbb{C}[P])
\end{align}
and $d$ preserves the grading. 
Define a degree for each monomial by 
$$\deg(x_\alpha)=1=\deg(\varphi^\alpha),\quad \deg(ab)=\deg(a)+\deg(b)$$
and denote by $F^pC^\bullet(\mathfrak{p};\mathbb{C}[P])$ the subspace spanned by all the monomials of degree greater than or equal to $p$.
Then \eqref{differential} implies $d\colon F^p C^\bullet(\mathfrak{p};\mathbb{C}[P])\rightarrow F^{p+1} C^\bullet(\mathfrak{p};\mathbb{C}[P])$.
Hence, there is a spectral sequence $E_r\Rightarrow H^\bullet(\mathfrak{p};\mathbb{C}[P])$ such that
$$E_1=H\left(\operatorname{gr}_FC^\bullet(\mathfrak{p};\mathbb{C}[P]); \operatorname{gr} d\right).$$ The filtration gives a filtration on each $C^\bullet_n(\mathfrak{p};\mathbb{C}[P])$ by 
$$F^pC^\bullet_n(\mathfrak{p};\mathbb{C}[P])=C^\bullet_n(\mathfrak{p};\mathbb{C}[P])\cap F^pC^\bullet(\mathfrak{p};\mathbb{C}[P]).$$
It is of finite length and satisfies $C^\bullet_n(\mathfrak{p};\mathbb{C}[P])=\cup_{p}F^pC^\bullet_n(\mathfrak{p};\mathbb{C}[P])$. Therefore, the spectral sequence $E_r$ converges. 
Since the differential on 
$\operatorname{gr}_FC^\bullet(\mathfrak{p};\mathbb{C}[P])$
is
$$\operatorname{gr} d=\sum_{\alpha=1}^N(-1)^{\bar{\alpha}}\frac{\partial}{\partial x_\alpha}\otimes \varphi^\alpha,$$
the complex $(\operatorname{gr}_FC^\bullet(\mathfrak{p};\mathbb{C}[P]),\operatorname{gr} d)$ is the (algebraic) de Rham complex of $P$. By Corollary \ref{unipotents are affine spaces}, 
$$H^n(\operatorname{gr}_FC^\bullet(\mathfrak{p};\mathbb{C}[P]),\operatorname{gr} d)\simeq H^n_{\mathrm{dR}}(\mathbb{C}^{p|q})$$
for some $p,q\in \mathbb{Z}_{\geq0}$. Then by Poincar\'{e} lemma below, we have $H^n_{\mathrm{dR}}(\mathbb{C}^{p|q})\simeq \delta_{n,0}\mathbb{C}$.
It follows that the spectral sequence collapses at $r=1$:
$$\operatorname{gr}_FH^n(\mathfrak{p};\mathbb{C}[P])\simeq E_1^n\simeq \delta_{n,0}\mathbb{C}.$$
This completes the proof.
\endproof
The following assertion is proved by Kostant in the analytic setting \cite[Theorem 4.6]{Ko} and is also well-known in the algebraic setting.
\begin{lemma}[Poincar\'{e} lemma]
$H^n_{\mathrm{dR}}(\mathbb{C}^{p|q})\simeq \delta_{n,0}\mathbb{C}$.
\end{lemma}
\proof 
The algebraic de Rham complex of $\mathbb{C}^{p|q}$ is 
$$C_{\mathrm {dR}}^\bullet (\mathbb{C}^{p|q})=\mathbb{C}[\mathbb{C}^{p|q}]\otimes S(\Pi \mathbb{C}^{p|q}).$$
We write $\mathbb{C}[\mathbb{C}^{p|q}]=\mathbb{C}[x_1,\cdots,x_p,\theta_1,\cdots,\theta_q]$ where $x_i$ (resp.\ $\theta_j$) are even (resp.\ odd) variables, and 
$S(\Pi \mathbb{C}^{p|q})=\mathbb{C}[dx_1,\cdots dx_p,d\theta_1,\cdots, d\theta_q]$
where $dx_i$ (resp.\ $d\theta_j$) are odd (resp.\ even) variables.
Then the differential is 
$$d=\sum_{i=1}^p \frac{\partial}{\partial x_i}\otimes dx_i-\sum_{j=1}^q \frac{\partial}{\partial \theta_j}\otimes d\theta_j.$$
Thus by K$\ddot{\text{u}}$nneth formula, 
$H^\bullet_{\mathrm{dR}}(\mathbb{C}^{p|q})\simeq \left(H^\bullet_{\mathrm{dR}}(\mathbb{C}^{1|0})\right)^{\otimes p}\otimes \left(H^\bullet_{\mathrm{dR}}(\mathbb{C}^{0|1})\right)^{\otimes q}$. Therefore, it suffices to show the assertion in the cases $\mathbb{C}^{1|0}$ and $\mathbb{C}^{0|1}$, respectively. The first one is the usual Poincar\'{e} lemma. For the second one, note that the $d$-closed forms are linear combinations of 
$1\otimes d\theta^n$, $(n\geq0)$, which are exact if and only if $n>0$ since
$$d(-\theta\otimes d\theta^{p-1})=1\otimes d\theta^p.$$
This completes the proof.
\endproof

\subsection{Superschemes of formal arcs}\label{Superscheme of formal arcs}
Let $\underbar{SSch}$ denote the category of superschemes over $\mathbb{C}$. An object 
$$D:=\mathrm{Spec}(\mathbb{C}[\![t]\!]),$$
of $\underbar{SSch}$ is called the formal disc.
\begin{proposition}[{\cite[Proposition 4.2.1]{KV}}]
Let $X$ be a superscheme over $\mathbb{C}$. The contravariant functor 
$$\underbar{SSch}\rightarrow \underbar{Set},\quad Y\mapsto \mathrm{Hom}_{\underbar{\tiny{SSch}}}(Y\widehat{\times} D,X)$$
is represented by a superscheme $JX$, that is,
$$\mathrm{Hom}_{\underbar{\tiny{SSch}}}(Y\widehat{\times} D,X)\simeq \mathrm{Hom}_{\underbar{\tiny{SSch}}}(Y,JX)$$
for any object $Y$ of $\underbar{SSch}$. Here $Y\widehat{\times} D$ is the completion of $Y\times D$ with respect to the subsuperscheme $Y\times\{0\}$.
\end{proposition}
The superscheme $JX$ in the above proposition is called the {\it superscheme of formal arcs} in $X$ or the {\it arc space} of $X$.
The association $X\mapsto JX$ gives a functor:
\begin{align}\label{arc space functor}
J\text{-}\colon \underbar{SSch}\rightarrow \underbar{SSch}.
\end{align}
The arc space $JX$ of $X$ has a canonical projection 
$\pi_X\colon JX \rightarrow X$ satisfying the following functoriality:
for any morphism $f\colon X\rightarrow Y$ of supershemes, the morphism $Jf\colon JX\rightarrow JY$ makes the following diagram commutative:
\begin{equation*}
\vcenter{\xymatrix{
JX \ar[r]^{Jf}\ar[d]^{\pi_X} &JY \ar[d]^{\pi_Y}\\
X\ar[r]^f &Y.}}
\end{equation*}
If $X$ is an affine superscheme with
$$\mathbb{C}[X]=\mathbb{C}[x_1,\cdots,x_N]/(f_1,\cdots,f_M),$$
where $\mathbb{C}[x_1,\cdots,x_N]$ is the ring of superpolynomials in the variables $x_1,\cdots,x_N$ and $f_1,\cdots,f_M\in \mathbb{C}[x_1,\cdots,x_N]$, then $JX$ is again an affine superscheme with 
$$\mathbb{C}[X]=\mathbb{C}[x_{1(n)},\cdots,x_{N{(n)}}\mid n<0]/( f_{1(n)},\cdots,f_{M(n)}\mid n<0).$$
$$\sum_{n<0}f_{j(n)}z^{-n-1}=f_j(x_1(z),\cdots,x_N(z)),$$
with $x_i(z)=\sum_{n<0}x_{i(n)}z^{-n-1}$. 
In this case, the canonical projection $\pi_X$ is given by
$$\mathbb{C}[X]\rightarrow \mathbb{C}[JX],\quad x_i\mapsto x_{i(-1)}.$$
Later, we use the following properties of the functor \eqref{arc space functor}.
\begin{lemma}[{cf. \cite{EM}}]\label{functoriality of arc}
For superschemes $X$ and $Y$, the following holds:
\begin{enumerate}
\item
$X\simeq Y$ implies $JX\simeq JY$.
\item
$J(X\times Y)\simeq JX\times JY$.
\end{enumerate}
\end{lemma}
For an algebraic supergroup $P$ with Lie superalgebra $\mathfrak{p}=T_eP$, the arc space $JP$ of $P$ has a structure of group superscheme and the Lie superalgebra on the tangent space $T_eJP$ is   $J \mathfrak{p}:=\mathfrak{p}[\![t]\!]$. The following proposition can be proved in the same way as Proposition \ref{cohomology of gra unipotent group}.
\begin{proposition}\label{cohomology vanishing for prounipotent AG}
For a positively graded unipotent affine algebraic supergroup $P$, we have
$$H^n(J \mathfrak{p};\mathbb{C}[JP])\simeq \delta_{n,0} \mathbb{C}.$$
\end{proposition}

\section{$\mathcal{W}$-superalgebras and Miura maps}\label{sec: W-superalgebra and Miura map} 

\subsection{$\mathcal{W}$-superalgebras}\label{W-superalgebras}
We review here the definition of (affine) $\mathcal{W}$-superalgebras \cite{FF0,KRW}. 
Following \cite{Ka3}, given a vertex superalgebra $V$, we denote by $\partial$ the translation operator and $:ab:$ the normally ordered product of $a,b\in V$, by $[a_\lambda b]=\sum_{n\geq0}  a_{(n)}b \lambda^n/n!$ the $\lambda$-bracket. We remark that if we denote by $a(z)=\sum_{n\in \mathbb{Z}}a_{(n)}z^{-n-1}$ (resp.\ $b(z)$) the field corresponding to $a$ (resp. $b$), then $:ab:=a_{(-1)}b$ and $[a_\lambda b]=\int e^{z \lambda} a(z)b\hspace{0.5mm}dz$.

Let $\mathfrak{g}$ be a simple finite-dimensional Lie superalgebra as in the introduction and $\mathfrak{g}=\mathfrak{g}_{\bar{0}}\oplus \mathfrak{g}_{\bar{1}}$ the parity decomposition.
Recall that the Lie superalgebra $\mathfrak{g}$ admits a non-degenerate even invariant bilinear form, which we denote by $\kappa$. Since such forms are all proportional, we may write $\kappa=k \kappa_0$ for some $k\in \mathbb{C}$ where $\kappa_0=(\cdot,\cdot)$ is the one satisfying $\kappa_0(\theta,\theta)=2$ for the highest root $\theta$ of $\mathfrak{g}_{\bar{0}}$, see \cite{Ka1, Mu}.
Let $f\in \mathfrak{g}_{\bar{0}}$ be a non-zero even nilpotent element. A good grading of $\mathfrak{g}$ with respect to $f$ is a $\frac{1}{2}\mathbb{Z}$-grading on $\mathfrak{g}$
$$\Gamma\colon \mathfrak{g}=\bigoplus_{j\in\frac{1}{2}\mathbb{Z}}\mathfrak{g}_j$$
such that $f\in \mathfrak{g}_{-1}$ and the adjoint action $\mathrm{ad}_f$ of $f$ is injective $\mathfrak{g}_j\hookrightarrow \mathfrak{g}_{j-1}$ for $j\geq1/2$ and surjective $\mathfrak{g}_j\twoheadrightarrow \mathfrak{g}_{j-1}$ for $j\leq1/2$, (see \cite{EK,H} for the definition and classification). We fix a basis $x_\alpha$, ($\alpha\in I=\{1,\cdots,\mathrm{dim}\mathfrak{g}\}$) of $\mathfrak{g}$ such that each $x_i$ is homogeneous with respect to the parity $\mathbb{Z}_2$ and the good grading $\Gamma$. Then we have $I=\sqcup_j I_j$ where $I_j=\{\alpha\mid x_\alpha\in \mathfrak{g}_j\}$. Let $c_{\alpha \beta}^\gamma$ denote the structure constants, i.e., $[x_\alpha,x_\beta]=\sum_{\gamma} c_{\alpha \beta}^\gamma x_\gamma$.

Let $V^k(\mathfrak{g})$ denote the universal affine vertex superalgebra of $\mathfrak{g}$ at level $k$, which is generated by $X$, $(X\in \mathfrak{g})$, satisfying the $\lambda$-bracket
$$[X_\lambda Y]=[X,Y] +k(X,Y)\lambda,\quad X,Y\in \mathfrak{g}.$$
We define a conformal grading on $V^k(\mathfrak{g})$ by $\operatorname{Delta}(u)=1-j$ for $u\in \mathfrak{g}_j$.
Let $F_{\mathrm{ch}}(\mathfrak{g}_+)$ be the charged fermion vertex superalgebra associated with $\mathfrak{g}_+:=\oplus_{j>0}\mathfrak{g}_j$. It is generated by $\varphi_{\alpha},\varphi^{\alpha}$, ($\alpha\in I_+:=\sqcup_{j>0}I_j$) of parity reversed to $x_\alpha$, satisfying the $\lambda$-bracket
\begin{align*}
[\varphi_{\alpha\lambda}\varphi^{\beta}]=\delta_{\alpha,\beta},\quad
[\varphi_{\alpha\lambda}\varphi_{\beta}]=0=[\varphi^{\alpha}\!_\lambda\varphi^{\beta}],\quad \alpha,\beta\in I_+.
\end{align*}
We define a conformal grading on $F_{\mathrm{ch}}(\mathfrak{g}_+)$ by $\operatorname{Delta}(\varphi_\alpha)=1-j$, $\operatorname{Delta}(\varphi^\alpha)=j$ for $\alpha\in I_j$, and 
a degree on $F_{\mathrm{ch}}(\mathfrak{g}_+)=\bigoplus_{n\in \mathbb{Z}}F_{\mathrm{ch}}^n$ by 
$$\deg(\varphi^{\alpha})=1=-\deg(\varphi_{\alpha}),\quad \alpha\in I_+,$$
$$\deg(\partial a)=\deg(a),\quad \deg(:ab:)=\deg(a)+\deg(b).$$
Let $\Phi(\mathfrak{g}_{1/2})$ be the neutral Fermion vertex superalgebra associated with $\mathfrak{g}_{1/2}$, which is generated by $\Phi_{\alpha}$, ($\alpha\in I_{1/2}$) satisfying the $\lambda$-bracket
\begin{align*}
[\Phi_{\alpha\hspace{0.5mm} \lambda}\Phi_{\beta}]=\chi([x_{\alpha},x_{\beta}]),\quad \alpha,\beta\in I_{\frac{1}{2}},
\end{align*}
where $\chi(x)=(f,x)$ for $x\in \mathfrak{g}$. 
We define a conformal grading on $\Phi(\mathfrak{g}_{1/2})$ by $\operatorname{Delta}(\Phi_\alpha)=1/2$ for all $\alpha\in I_{1/2}$.
Introduce a $\mathbb{Z}$-graded vertex superalgebra by 
\begin{align*}
C_k^\bullet(\mathfrak{g},f;\Gamma)=V^k(\mathfrak{g}) \otimes \Phi(\mathfrak{g}_{\frac{1}{2}}) \otimes F_\mathrm{ch}^\bullet(\mathfrak{g}_+),
\end{align*}
with cohomological grading given by the degree on $F_\mathrm{ch}(\mathfrak{g}_+)$ and define a differential on $C_k^\bullet(\mathfrak{g},f;\Gamma)$ by $d_{(0)}$ of the element
\begin{align*}
d=\sum_{\alpha\in I_+}:((-1)^{\bar{\alpha}} x_\alpha+&\Phi_\alpha+\chi(x_\alpha))\varphi^{\alpha}: -\frac{1}{2}\sum_{\alpha,\beta,\gamma\in I_+}(-1)^{\bar{\alpha}\bar{\gamma}}c_{\alpha\beta}^{\gamma}:\varphi_{\gamma}\varphi^{\alpha}\varphi^{\beta}:,
\end{align*}
where $\Phi_\alpha=0$ for $\alpha\notin I_{1/2}$. Then $(C_k^\bullet(\mathfrak{g},f;\Gamma),d_{(0)})$ forms a cochain complex, called the {\it BRST complex}. The vertex superalgebra obtained as the cohomology $H(C_k^\bullet(\mathfrak{g},f;\Gamma),d_{(0)})$ is called the (affine) {\it $\mathcal{W}$-superalgebra} associated with $(\mathfrak{g},f,k,\Gamma)$ and denoted by $\mathcal{W}^k(\mathfrak{g},f;\Gamma)$. By \cite[Theorem 4.1]{KW1}, we have
$$\mathcal{W}^k(\mathfrak{g},f;\Gamma)=H^0(C_k^\bullet(\mathfrak{g},f;\Gamma),d_{(0)}).$$
Since $d_{(0)}$ preserves the conformal grading on $C_k^\bullet(\mathfrak{g},f;\Gamma)$ induced by those on each component $V^k(\mathfrak{g})$, $\Phi(\mathfrak{g}_{1/2})$ and $F_\mathrm{ch}^\bullet(\mathfrak{g}_+)$, $\mathcal{W}^k(\mathfrak{g},f;\Gamma)$ has an induced conformal grading, which is a $\frac{1}{2}\mathbb{Z}_{\geq0}$-grading.
 
\subsection{Miura maps}
Define an element
$$J^u=u+\sum_{\alpha,\beta\in I_+}(-1)^{\bar{\alpha}}c_{u,\beta}^\alpha:\varphi_\alpha\varphi^\beta:,\quad u\in \mathfrak{g}.$$
Let $C_+\subset C_k^\bullet(\mathfrak{g},f;\Gamma)$ denote the vertex subalgebra generated by $J^u$, ($u\in \mathfrak{g}_+$), and $\varphi_\alpha$, ($\alpha\in I_+$), and let $C_-\subset C_k^\bullet(\mathfrak{g},f;\Gamma)$ denote the one generated by $J^u$, ($u\in \mathfrak{g}_{\leq0}:=\oplus_{j\leq0}\mathfrak{g}_j$), $\Phi_\alpha$, ($\alpha\in I_{1/2}$), and $\varphi^\alpha$, ($\alpha\in I_+$).
By \cite{KW1,KW2}, $C_\pm^\bullet\subset C_k^\bullet(\mathfrak{g},f;\Gamma)$ are subcomplexes and give a decomposition 
$$C_k^\bullet(\mathfrak{g},f;\Gamma)\simeq C_+^\bullet\otimes C_-^\bullet.$$
Moreover, we have $H^n(C_+^\bullet,d_{(0)})\simeq \delta_{n,0} \mathbb{C}$. Thus we have an isomorphism
$$\mathcal{W}^k(\mathfrak{g},f;\Gamma)\simeq H^0(C_-^\bullet,d_{(0)}).$$
Since $C_-$ is $\mathbb{Z}_{\geq0}$-graded as a complex, it follows that $\mathcal{W}^k(\mathfrak{g},f;\Gamma)$ is a vertex subalgebra of $C_-^0$. By \cite[Theorem 2.4, (c)]{KRW}, we have
$$[J^u_\lambda J^v]=J^{[u,v]}+\tau_k(u,v)\lambda,\quad u,v\in \mathfrak{g}_{\leq0},$$
where
$$\tau_k(u,v)=k(u,v)+\frac{1}{2}\left(\kappa_\mathfrak{g}(u,v)-\kappa_{\mathfrak{g}_0}(u,v)\right)$$
and $\kappa_\mathfrak{g}$ (resp.\ $\kappa_{\mathfrak{g}_0}$) denotes the Killing form of $\mathfrak{g}$ (resp.\ $\mathfrak{g}_0$). 
Therefore, the elements $J^u$, $u\in \mathfrak{g}_{\leq0}$, (resp.\ $u\in \mathfrak{g}_{<0}$, $u\in \mathfrak{g}_{0}$) generate the universal affine vertex superalgebra $V^{\tau_k}(\mathfrak{g}_{\leq0})$ (resp.\ $V^{\tau_k}(\mathfrak{g}_{<0})$, $V^{\tau_k}(\mathfrak{g}_{0})$) of $\mathfrak{g}_{\leq0}$ (resp.\ $\mathfrak{g}_{<0}$, $\mathfrak{g}_0$) at level $\tau_k$.
It follows that
$$C_-^0\simeq V^{\tau_k}(\mathfrak{g}_{\leq0})\otimes \Phi(\mathfrak{g}_{\frac{1}{2}}).$$
Since $V^{\tau_k}(\mathfrak{g}_{<0})\otimes \Phi(\mathfrak{g}_{1/2})\subset V^{\tau_k}(\mathfrak{g}_{\leq0})\otimes \Phi(\mathfrak{g}_{\frac{1}{2}})$ is a (vertex superalgebra) ideal, we have a natural surjection
$$V^{\tau_k}(\mathfrak{g}_{\leq0})\otimes \Phi(\mathfrak{g}_{\frac{1}{2}})\twoheadrightarrow \left(V^{\tau_k}(\mathfrak{g}_{\leq0})\otimes \Phi(\mathfrak{g}_{\frac{1}{2}})\right)/\left(V^{\tau_k}(\mathfrak{g}_{<0})\otimes \Phi(\mathfrak{g}_{1/2})\right)\simeq V^{\tau_k}(\mathfrak{g}_0)\otimes \Phi(\mathfrak{g}_{\frac{1}{2}}).$$
The restriction to the subalgebra $\mathcal{W}^k(\mathfrak{g},f;\Gamma)$ is called the {\it Miura map} \cite{A,FBZ,FF2,G}
\begin{align}\label{eq: Miura map}
\mu\colon \mathcal{W}^k(\mathfrak{g},f;\Gamma)\rightarrow V^{\tau_k}(\mathfrak{g}_0)\otimes \Phi(\mathfrak{g}_{\frac{1}{2}}).
\end{align}
Note that $V^{\tau_k}(\mathfrak{g}_0)\otimes \Phi(\mathfrak{g}_{\frac{1}{2}})$ has a conformal grading defined by 
$$\Delta(J^u(z))=1,\ (u\in \mathfrak{g}_0),\quad \Delta(\Phi_\alpha(z))=\frac{1}{2},\ (\alpha\in I_{\frac{1}{2}}),$$
and that $\mu$ preserves the conformal grading $\Delta$. 
\begin{theorem}\label{th: injectivity of Miura}
The Miura map $\mu$ is injective.
\end{theorem}
The above theorem is proven in the literature under some conditions: \cite{F,A} provide a proof in the non-super setting with $k$ arbitrary and $f$ principal nilpotent, \cite[Corollary 8.2]{DKV} in the non-super setting in full generality, and \cite{G} in the super-setting with $k$ generic. 
The proof in \cite{A} is geometric and, indeed, applies for an arbitrary nilpotent element $f$ in the non-super setting.
We prove Theorem \ref{th: injectivity of Miura}, essentially following \cite{A} together with necessary supergeometry.

\subsection{Li's filtration}
By \cite{L}, given a vertex superalgebra $V$, the subspaces $F^pV$, ($p\in \mathbb{Z}$), spanned by
\begin{align*}
a^1_{(-n_1-1)}a^2_{(-n_2-1)}\cdots a^r_{(-n_r-1)}|0\rangle
\end{align*}
with $a^1, a^2, \cdots, a^r\in V$, $n_i \geq 0$, $n_1+n_2+\cdots +n_r\geq p$, form a descending filtration $F^\bullet V$ of $V$ satisfying the following propertires:
\begin{itemize}
\item
$F^pV_{(n)}F^qV\subset F^{p+q-n-1}V$, ($p,q\in \mathbb{Z}_{\geq0},\ n\in \mathbb{Z}$),
\item
$F^pV_{(n)}F^qV\subset F^{p+q-n}V$, ($p,q\in \mathbb{Z}_{\geq0},\ n\geq0$).
\end{itemize}
Let 
$$\operatorname{gr}_F V:=\oplus_{p\in \mathbb{Z}} \operatorname{gr}_F^p V,\quad \operatorname{gr}_F^pV=F^pV/F^{p+1}V,$$
denote the associated graded vector superspace and $\sigma_p\colon F^pV\twoheadrightarrow \operatorname{gr}_F^pV$ the canonical projection. 
By \cite{L}, $\operatorname{gr}_F V$ has a Poisson vertex superalgebra structure by 
\begin{align*}
&\partial \sigma_p(a)=\sigma_{p+1}(\partial a),\quad \sigma_p(a)\sigma_q(b)=\sigma_{p+q}(a_{(-1)}b),\\
&\{\sigma_p(a)_\lambda\sigma_q(b)\}=\sum_{n\geq0}\frac{1}{n!}\sigma_{p+q-n}(a_{(n)}b)\lambda^n.
\end{align*}
(see e.g.\ \cite{FBZ,Ka3} for the original definition of Poisson vertex superalgebras and, especially, \cite{Ka3} in terms of the $\lambda$-bracket). 
It is obvious that any homomorphism $\eta\colon V\rightarrow W$ of vertex superalgebras induces a homomorphism of Poisson vertex superalgebras
$$\bar{\eta}\colon  \operatorname{gr}_F V\rightarrow \operatorname{gr}_F W,\quad \sigma_p(a)\mapsto \sigma_p(\eta(a)).$$
\begin{lemma}\label{injectivity: general setting}
Let $V$ and $W$ be vertex superalgebras equipped with conformal $\frac{1}{2}\mathbb{Z}_{\geq0}$-gradings and $\eta\colon V\rightarrow W$ be a homomorphism of vertex superalgebras preserving the conformal gradings. Then $\eta$ is injective if $\bar{\eta}$ is injective.
\end{lemma}
\proof
Let $V=\oplus_{\Delta\in\frac{1}{2}\mathbb{Z}_{\geq0}}V_\Delta$ (resp.\ $W=\oplus_{\Delta\in\frac{1}{2}\mathbb{Z}_{\geq0}}W_\Delta$) denote the conformal grading on $V$ (resp\ $W$).
It is immediate from the definition that Li's filtration $F^\bullet V$ defines a filtration on each $V_\Delta$ by setting
$$F^pV_\Delta=V_\Delta\cap F^pV.$$
Moreover, it is of finite length, i.e., $V_\Delta\cap F^pV=0$ for $p\gg0$.
Since $\eta$ preserves the conformal gradings, $\operatorname{gr}_F \eta$ is restricted to 
$$\operatorname{gr}_F \eta\colon \operatorname{gr}_F V_\Delta\rightarrow \operatorname{gr}_F W_\Delta,\quad \Delta\geq0.$$
It is straightforward to show the injectivity of $\eta$ on $V_\Delta$ inductively from that of $\operatorname{gr}_F\eta$ on $\operatorname{gr}_F V_\Delta$ starting with the subspace $\operatorname{gr}_F^N V_\Delta$ such that $\operatorname{gr}_F^nV_\Delta=0$ for $n>N$.
\endproof
By \cite[Section 3]{A}, the Poisson vertex superalgebra $\operatorname{gr}_FV^{\tau_k}(\mathfrak{g}_0)\otimes \Phi(\mathfrak{g}_{1/2})$ is isomorphic to the algebra of differential superpolynomials 
$$S^\partial\left(\mathfrak{g}_0\oplus \mathfrak{g}_{\frac12}\right)=\mathbb{C}\left[\partial^n x_\alpha\left|\, \alpha\in I_0\sqcup I_{\frac12},n\geq0\right.\right]$$
with $\lambda$-bracket
\begin{align*}
&\{x_{\alpha\lambda} x_\beta\}=[x_\alpha,x_\beta],\quad \alpha,\beta\in I_0,\\
&\{x_{\alpha\lambda} x_{\beta}\}=\chi([x_\alpha,x_\beta]),\quad \alpha,\beta\in I_{\frac{1}{2}},\\
&\{x_{\alpha\lambda} x_\beta\}=0,\quad \alpha\in I_0,\ \beta\in I_{\frac{1}{2}}.
\end{align*}
Thus the Miura map $\mu$ in \eqref{eq: Miura map} induces a homomorphism of Poisson vertex superalgebras
\begin{align}\label{eq: graded Miura map}
\bar{\mu}\colon \mathrm{gr}_F\mathcal{W}^k(\mathfrak{g},f,\Gamma)\rightarrow S^\partial\left(\mathfrak{g}_0\oplus \mathfrak{g}_{\frac12}\right).
\end{align}
Lemma \ref{injectivity: general setting} implies the following.
\begin{corollary}\label{corollary: implication of injectivity}
The injectivity of $\bar{\mu}$ implies that of $\mu$.
\end{corollary}
\subsection{Arc spaces of Slodowy slices}\label{sec: Arc space of Slodowy slice}
Consider Li's filtration $F^\bullet C_-^\bullet$ on the complex $C_-^\bullet$. The associated graded vector superspace $\bar{C}_-^\bullet=\operatorname{gr}_FC_-^\bullet$ is again a complex with differential $\bar{d}_{(0)}$. Note that $\operatorname{gr}_FC_-^\bullet$ is algebra of differential superpolynomials
$$\operatorname{gr}_FC_-^\bullet\simeq S^\partial(\mathfrak{g}_{\leq0}\oplus \mathfrak{g}_{\frac{1}{2}})\otimes S^\partial(\Pi \mathfrak{g}_+^*),$$
equipped with $\lambda$-bracket
\begin{align*}
&\{u_\lambda v\}=[u,v],\quad (u,v\in \mathfrak{g}_{\leq0}),\\
&\{u_\lambda v\}=\chi([u,v]),\quad (u,v\in \mathfrak{g}_{\frac12}),\\
&\{\varphi^\alpha_\lambda u\}=\sum_{\beta\in I_+}c_{u,\beta}^\alpha\varphi^\alpha,\quad (\alpha\in I_+, u\in \mathfrak{g}_{\leq0}),\\
&\{u_\lambda v\}=0=\{\varphi^\alpha_\lambda\varphi^\beta\},\quad (u\in \mathfrak{g}_{\leq0},\ v\in \mathfrak{g}_{\frac12},\ \alpha,\beta\in I_+),
\end{align*}
where $\{\varphi^\alpha\}_{\alpha\in I_+}\subset \mathfrak{g}_+^*$ is the dual basis of $\{x_\alpha\}_{\alpha\in I_+}$ of $\mathfrak{g}_+$ with reversed parity.
The differential $Q=\bar{d}_{(0)}$ is given by
\begin{align}
&Q(u)=
\begin{cases}
\displaystyle{\sum_{\alpha\in I_+}\varphi^\alpha\left((-1)^{\bar{\alpha}}[x_\beta,u]_{\leq0}-(-1)^{\bar{u}}\left([x_\alpha,u]_{1/2}+\chi([x_\alpha,u]_{\geq1})\right)\right)},& (u\in \mathfrak{g}_{\leq0}),\\
\displaystyle{\sum_{\alpha\in I_+}\chi([x_\alpha,u])\varphi^\alpha},& (u\in \mathfrak{g}_{1/2}),
\end{cases}\\
&Q(\varphi^\alpha)=-\frac{1}{2}\sum_{\beta,\gamma\in I_+}(-1)^{\bar{\alpha}\bar{\gamma}}c_{\alpha,\beta}^\gamma \varphi^\beta \varphi^\gamma,\quad \alpha\in I_+,\\
&Q(\partial a)=\partial Q(a),\quad Q(ab)=Q(a)b+(-1)^{\bar{a}}aQ(b).
\end{align}
Here we have used the projections 
$$\mathfrak{g}\rightarrow \mathfrak{g}_{\leq0}\oplus \mathfrak{g}_{1/2}\oplus \mathfrak{g}_{\geq1},\quad u\mapsto (u_{\leq0},u_{1/2},u_{\geq1}).$$
We will interpret the complex $(\operatorname{gr}_FC_-^\bullet,Q)$ geometrically. To this end, let us first consider the finite analogue: the quotient graded superspace
$$\operatorname{gr}_F^\mathrm{fin}C_-^\bullet:=\operatorname{gr}_FC_-^\bullet/(a\partial(b)\mid a,b\in \operatorname{gr}_FC_-^\bullet).$$
It is easy to check that for a vector superspace $A$, 
$$S^\partial(A)/(a\partial(b)\mid a,b\in A)\simeq S(A),$$
where $S(A)$ is the symmetric superalgebra of $A$. Thus we have
$$\operatorname{gr}_F^\mathrm{fin}C_-\simeq S(\mathfrak{g}_{\leq0}\oplus \mathfrak{g}_{1/2})\otimes S(\Pi \mathfrak{g}_+^*).$$
Let $\operatorname{gr}_FC_-\rightarrow \operatorname{gr}_F^\mathrm{fin}C_-$, ($a\mapsto [a]$), denote the canonical projection. It is naturally a Poisson superalgebra by
$$[a][b]=[ab],\quad \left\{[a],[b]\right\}=[\{a_\lambda b\}|_{\lambda=0}].$$
The Poisson superalgebra $\operatorname{gr}_F^{\mathrm{fin}}C_-^\bullet$ is called the {\em Zhu's $C_2$-algebra} of $\operatorname{gr}_FC_-^\bullet$. By (2.5), the differential $Q$ induces a differential $Q^{\mathrm{fin}}$ on $\operatorname{gr}_F^{\mathrm{fin}}C_-^\bullet$, which is determined by the formulas (2.3)-(2.5) with $u$ replaced by $[u]$, e.g.,
$$Q^{\mathrm{fin}}([\varphi^\alpha])=-\frac{1}{2}\sum_{\beta,\gamma\in I_+}(-1)^{\bar{\alpha}\bar{\gamma}}c_{\alpha,\beta}^\gamma [\varphi^\beta] [\varphi^\gamma],\quad \alpha\in I_+.$$
Then the complex $(\operatorname{gr}_F^{\mathrm{fin}}C_-^\bullet,Q^{\mathrm{fin}})$ is identified with the Chevalley-Eilenberg complex of $\mathfrak{g}_+$ with coefficients in $S(\mathfrak{g}_{\leq0}\oplus \mathfrak{g}_{1/2})$ since the differential $Q^{\mathrm{fin}}$ defines a left $\mathfrak{g}_+$-module structure on $S(\mathfrak{g}_{\geq0}\oplus \mathfrak{g}_{1/2})$ by
\begin{align*}
x_\alpha\cdot u=
\begin{cases}
(-1)^{\bar{\alpha}}[x_\beta,u]_{\leq0}-(-1)^{\bar{u}}\left([x_\alpha,u]_{1/2}+\chi([x_\alpha,u]_{\geq1})\right),& (u\in \mathfrak{g}_{\leq0}),\\
\chi([x_\alpha,u]),& (u\in \mathfrak{g}_{1/2}),
\end{cases}
\end{align*}
for $\alpha\in I_+$. Let $G_+$ be the unipotent algebraic supergroup whose Lie superalgebra is $\mathfrak{g}_+$. The right $G_+$-action 
$$(f+\mathfrak{g}_{\geq-1/2})\times G_+\rightarrow (f+\mathfrak{g}_{\geq-1/2}),\quad (f+X,g)\mapsto g^{-1}(f+X)g,$$
induces a left $\mathfrak{g}_+$-module structure on the coordinate ring $\mathbb{C}[f+\mathfrak{g}_{\geq-1/2}]$. Then the isomorphism of $\mathbb{C}$-superalgebras
\begin{align*}
S(\mathfrak{g}_{\leq0}\oplus&\mathfrak{g}_{1/2})\simeq \mathbb{C}[f+\mathfrak{g}_{\geq-1/2}],\\
&u\mapsto 
\begin{cases}
-(-1)^{\bar{u}}(u|?),& u\in \mathfrak{g}_{\leq0},\\
(u|?),& u\in \mathfrak{g}_{1/2},
\end{cases}
\end{align*}
is a $\mathfrak{g}_+$-homomorphism.  Thus we have proved the following.
\begin{lemma}\label{identification of f.d. complex}
The complex $(\operatorname{gr}_F^{\mathrm{fin}}C_-^\bullet,Q^{\mathrm{fin}})$ is quasi-isomorphic to the Chevalley-Eilenberg complex of $\mathfrak{g}_+$ with coefficients in $\mathbb{C}[f+\mathfrak{g}_{\geq-1/2}]$.
\end{lemma}
By Jacobson-Morozov theorem, there exits an $\mathfrak{sl}_2$-triple $\{e,h,f\}\subset \mathfrak{g}_{\bar{0}}$ containing $f$. Set $\mathfrak{g}^e=\{X\in \mathfrak{g}\mid [e,X]=0\}$. The subvariety $\mathcal{S}_f=f+\mathfrak{g}^e\subset \mathfrak{g}$ is called the {\it Slodowy slice} of $\mathfrak{g}$ associated with $f$.
\begin{proposition}\label{slodowy slice}
We have an isomorphism
\begin{align}\label{eq: Gan-Ginzburg}
\xi\colon \mathcal{S}_f\times G_+\simeq f+\mathfrak{g}_{\geq-\frac{1}{2}},\quad (X,g)\mapsto g^{-1}Xg
\end{align}
of affine supervarieties.
\end{proposition}
\proof
It suffices to show that $\xi$ gives an isomorphism of all the $A$-valued points for an arbitrary commutative $\mathbb{C}$-superalgebra $A$. 
Let $f+X\in f+\mathfrak{g}^e(A)$ and $g\in G_+(A)$. Let $Y\in \mathfrak{g}_+(A)$ denote the element corresponding to $g$ by the isomorphism $\mathfrak{g}_+\simeq G_+$ as affine supervarieties by Corollary \ref{unipotents are affine spaces}. It satisfies
$$Z:=g^{-1}(f+X)g=\sum_{n\geq0}\frac{1}{n!} [\cdots [[f+X,\underbrace{Y],Y]\cdots Y]}_{n}.$$
For $u\in \mathfrak{g}$, decompose it as $u=\sum_pu_p$ with respect to the grading $\Gamma$. Then we have
\begin{align*}
Z_p=
\begin{cases}
0,&(p<-1),\\
f,&(p=-1),\\
X_p+[f,Y_{p+1}]+R_p,&(p\geq-\frac{1}{2}),
\end{cases}
\end{align*}
where $R_p$ is the term determined by $X_i$, ($i<p$), and $Y_j$, ($j<p+1$). Now, the assertion follows from the decomposition
$$\mathfrak{g}_p=\mathfrak{g}_p^e\oplus [f,\mathfrak{g}_{p+1}],\quad p\geq-1/2,$$
which is an immediate consequence of the definition of the good grading $\Gamma$.
\endproof
Note that the Poisson superalgeba structure on $S(\mathfrak{g}_{\leq0}\oplus\mathfrak{g}_{1/2})\simeq \mathbb{C}[f+\mathfrak{g}_{\geq-1/2}]	$ restricts to $\mathbb{C}[\mathcal{S}_f]\simeq \mathbb{C}[f+\mathfrak{g}_{\geq-1/2}]^{G_+}$. Thus $\mathcal{S}_f$ is a Poisson supervariety. We remark that if $\Gamma$ is a $\mathbb{Z}$-grading, then $\mathcal{S}_f$ is the Poisson supervariety obtained as the Hamiltonian reduction of $\mathfrak{g}$ with respect to the adjoint $G_+$-action and an infinitesimal character $\chi=(f|\text{-}):\mathfrak{g}_+\rightarrow \mathbb{C}$.
\begin{corollary}\label{corollary: arc space of Slodowy slice}\hspace{1mm}\\
(1) $H^n(\operatorname{gr}_F^{\mathrm{fin}}C_-^\bullet,Q^{\mathrm{fin}})=0$ for $n\neq0$.\\
(2) $H^0(\operatorname{gr}_F^{\mathrm{fin}}C_-^\bullet,Q^{\mathrm{fin}})\simeq \mathbb{C}[\mathcal{S}_f]$ as Poisson superalgebras.
\end{corollary}
\proof
By Lemma \ref{identification of f.d. complex} and Proposition \ref{slodowy slice}, we have
$$H^n(\operatorname{gr}_F^{\mathrm{fin}}C_-^\bullet,Q^{\mathrm{fin}})\simeq H^n(\mathfrak{g}_+;\mathbb{C}[f+\mathfrak{g}_{\geq-1/2}])\simeq H^n(\mathfrak{g}_+;\mathbb{C}[G_+])\otimes \mathbb{C}[\mathcal{S}_f].$$
Since $G_+$ is a positively graded unipotent algebraic supergroup, the assertion follows from Proposition \ref{cohomology of gra unipotent group}.
\endproof

Next, we consider the complex $(\operatorname{gr}_F C_-^\bullet,Q)$. Recall that the arc space $JG_+$ of $G_+$ is a group superscheme whose Lie superalgebra is $J \mathfrak{g}_+=\mathfrak{g}_+[\![t]\!]$. 
It follows from Lemma \ref{identification of f.d. complex} that $(\operatorname{gr}_FC^\bullet_-,Q)$ is quasi-isomorphic to the Chevalley-Eilenberg complex of $J \mathfrak{g}_+$ with coefficients in $\mathbb{C}[f+J \mathfrak{g}_{\geq-1/2}]$ with left $J \mathfrak{g}_+$-module structure induced by the right $JG_+$-action
\begin{align}\label{jet action}
J(f+\mathfrak{g}_{\geq-\frac{1}{2}})\times JG_+\rightarrow J(f+\mathfrak{g}_{\geq-\frac{1}{2}}),\quad (f+X,g)\mapsto g^{-1}(f+X)g.
\end{align}
The following is a generalization of \cite[Theorem 5.7]{A}.
\begin{lemma}\label{graded cohomology}\hspace{1mm}\\
(1) $H^n(\operatorname{gr}_FC^\bullet_-,Q)=0$ for $n\neq0$.\\
(2) $H^0(\operatorname{gr}_FC^\bullet_-,Q)\simeq \mathbb{C}[J\mathcal{S}_f]$ as Poisson vertex superalgebras.
\end{lemma}
\proof
By Lemma \ref{functoriality of arc} and Proposition \ref{slodowy slice}, we have an isomorphism
$$J\mathcal{S}_f\times JG_+\simeq J(f+\mathfrak{g}_{\geq-1/2}).$$
Therefore, by using Proposition \ref{cohomology vanishing for prounipotent AG}, we have 
\begin{align*}
H^n(\operatorname{gr}_FC^\bullet_-,A)
&\simeq H^n(J \mathfrak{g}_+;\mathbb{C}[f+J \mathfrak{g}_{\geq-1/2}])\\
&\simeq H^n(J \mathfrak{g}_+;\mathbb{C}[JG_+])\otimes \mathbb{C}[J\mathcal{S}_f]\\
&\simeq \delta_{n,0} \mathbb{C}[J\mathcal{S}_f].
\end{align*}
\endproof
The following is a generalization of \cite[Theorem 5.8]{A}.
\begin{proposition}\label{proposition: classical W-superalgebra}
$$\mathrm{gr}_F\mathcal{W}^k(\mathfrak{g},f,\Gamma)\simeq \mathbb{C}[J\mathcal{S}_f]$$
as Poisson vertex superalgebras.
\end{proposition}
\proof
The filtration $F^\bullet C_-^\bullet$ induces a filtration on $H\left(C_-^\bullet,d\right)$ by 
$$\mathcal{F}^nH\left(C_-^\bullet,d\right)=\mathrm{Im}\left(\mathrm{Ker}(d)\cap F^nC_-^\bullet\rightarrow H(C_-^\bullet,d)\right)$$
and thus a spectral sequence $E_1=H(\operatorname{gr}_FC^\bullet_-,Q)\Rightarrow H(C_-^\bullet,Q)$.
Recall that $C_-^\bullet$ is $\frac12 \mathbb{Z}_{\geq0}$-graded $C_-^\bullet=\oplus_{\Delta\in\frac{1}{2}\mathbb{Z}_{\geq0}}C_{-,\Delta}^\bullet$ and that $F^\bullet C_-^\bullet$ induces a filtration on each $C_{-,\Delta}^\bullet$ by 
$$F^pC_{-,\Delta}^\bullet=F^pC_{-}^\bullet\cap C_{-,\Delta}^\bullet.$$
Since it is of finite length, the spectral sequence $E_r$ converges.
By Lemma \ref{graded cohomology}, $E_r$ collapses at $r=1$:
$$\mathrm{gr}_\mathcal{F}H^\bullet(C_-^\bullet,d)\simeq E_1^\bullet=H^\bullet(\operatorname{gr}_FC^\bullet_-,Q)\simeq \delta_{n,0}\mathbb{C}[\mathcal{S}_f].$$
Therefore, it suffices to show 
$\mathrm{gr}_F\mathcal{W}^k(\mathfrak{g},f;\Gamma)\simeq \mathrm{gr}_\mathcal{F}H^0(C_-^\bullet,d)$.
Since we have a map
$F^n\mathcal{W}^k(\mathfrak{g},f;\Gamma)\rightarrow \mathcal{F}^nH^0(C_-^\bullet,d)$
by construction, we obtain a homomorphism of Poisson vertex superalgebras 
$$\mathrm{gr}_F\mathcal{W}^k(\mathfrak{g},f;\Gamma)\rightarrow \mathrm{gr}_\mathcal{F}H^0(C_-^\bullet,d)\simeq \mathbb{C}[J\mathcal{S}_f].$$
By \cite[Theorem 4.1]{KW1}, $\mathcal{W}^k(\mathfrak{g},f;\Gamma)$ is strongly generated by a basis of $\mathfrak{g}^e\simeq f+\mathfrak{g}^e=\mathcal{S}_f$ and we have a surjection of Poisson vertex superalgebras
$$\mathbb{C}[J\mathcal{S}_f]\twoheadrightarrow \mathrm{gr}_F\mathcal{W}^k(\mathfrak{g},f;\Gamma),$$
by \cite[Theorem 10]{A}. Since the composition of the above two homomorphisms is the identity of $\mathbb{C}[J\mathcal{S}_f]$, we obtain the assertion.
\endproof
\subsection{Injectivity of Miura maps}
It follows from Proposition \ref{proposition: classical W-superalgebra} and the proof of Corollary \ref{corollary: arc space of Slodowy slice} that the map $\bar{\mu}$ in \eqref{eq: graded Miura map} is identified with the composition 
\begin{align*}
\mathbb{C}[J\mathcal{S}_f]\hookrightarrow \mathbb{C}[J\mathcal{S}_f]\otimes \mathbb{C}[JG_+]
\simeq \mathbb{C}[J(f+\mathfrak{g}_{\geq-1/2})]\twoheadrightarrow \mathbb{C}[J(f+\mathfrak{g}_{\mathrm{ini}})],
\end{align*}
where $\mathfrak{g}_{\mathrm{ini}}=\mathfrak{g}_{-1/2}\oplus \mathfrak{g}_{0}$.
Note that it is the image of the functor $J\text{-}$ in \eqref{arc space functor} of its finite analogue
\begin{align}\label{eq: finite analogue of graded Miura}
\bar{\mu}_{\text{fin}}\colon \mathbb{C}[\mathcal{S}_f]\hookrightarrow \mathbb{C}[\mathcal{S}_f]\otimes \mathbb{C}[G_+]
\simeq \mathbb{C}[f+\mathfrak{g}_{\geq-\frac{1}{2}}]\twoheadrightarrow \mathbb{C}[f+\mathfrak{g}_{\mathrm{ini}}].
\end{align}
Since the injectivity of $\bar{\mu}_{\text{fin}}$ implies that of $\bar{\mu}$, the proof of Theorem \ref{th: injectivity of Miura} is reduced to the following by Corollary \ref{corollary: implication of injectivity}.
\begin{proposition}\label{classical finite Miura}
The map $\bar{\mu}_{\operatorname{fin}}$ is injective.
\end{proposition}
\proof
Recall that we have a functor from the category of affine supervarieties to the category of affine varieties of taking the reduced variety:
$$X=\mathrm{Spec}(\mathbb{C}[X])\mapsto X_{\mathrm{red}}:=\mathrm{Spec}\left(\mathbb{C}[X]/(\mathbb{C}[X]_{\bar{1}})\right)$$
where $(\mathbb{C}[X]_{\bar{1}})\subset \mathbb{C}[X]$ denotes the ideal generated by the odd subspace of $\mathbb{C}[X]$, see \cite{CCF}. Applying this functor to $\bar{\mu}_{\text{fin}}$, we obtain the Miura map for the reductive Lie algebra $\mathfrak{g}_{\bar{0}}$ with the same nilpotent element $f$. By \cite[Lemma 5.12]{G2}, it is injective. Note that the map $\bar{\mu}_{\text{fin}}$ preserves the parity. Then it suffices to show that the images of the linear coordinates of the odd part $\mathcal{S}_{f,\bar{1}}$ are linearly independent.
Therefore, it suffices to show that the image of the map 
\begin{align}\label{eq: needed property}
G_{+,\bar{0}}\times (f+\mathfrak{g}_{\operatorname{ini},\bar{1}})\rightarrow f+\mathfrak{g}_{\geq-1/2}\twoheadrightarrow \mathcal{S}_{f,\bar{1}}
\end{align}
is dense in the Zariski topology. Here the last projection is the composition
\begin{align*}
\mathfrak{g}_{\geq-1/2}
&=\mathfrak{g}_{\geq-1/2,\bar0}\oplus \mathfrak{g}_{\geq-1/2,\bar1}\\
&\twoheadrightarrow \mathfrak{g}_{\geq-1/2,\bar1}=
\mathfrak{g}_{\geq-1/2,\bar1}^e\oplus \left(\mathfrak{g}_{\geq-1/2,\bar1}\cap \mathrm{Im}(\mathrm{ad}_f)\right)\twoheadrightarrow \mathfrak{g}_{\geq-1/2,\bar1}^e.
\end{align*}
Since the restriction of \eqref{eq: needed property} to 
$$\{\mathrm{exp}(e)\}\times (f+\mathfrak{g}_{\operatorname{ini},\bar{1}})\rightarrow \mathcal{S}_{f,\bar{1}}$$
is an isomorphism as affine varieties, we obtain the assertion.
\endproof


\end{document}